\documentclass[10pt]{article}
\usepackage{graphicx}
\usepackage{todonotes}

\def\R{{\rm I\! R}}

\def\Z{{\bf \rm Z}}

\def \pTW*{\partial_{W^*} T}
\def \pTnW*{\partial^{(n)}_{W^*} T}

\newtheorem{theorem}{Theorem}

\title{Every period annulus is both reversible and symmetric } 

\author{Marco Sabatini 
\footnote{Dip. di Matematica, Univ. di Trento, I-38123 Povo, (TN) - Italy.
Phone: ++39(0461)281670, Fax: ++39(0461)281624, Email: marco.sabatini@unitn.it - \ \ \ \ \ \ \ \  \ \ \ \ \ \ \ \ \ \ \ \ This paper was partially supported by the GNAMPA, {\it Gruppo Nazionale per l'Analisi Matematica, la Probabilit\`a e le loro Applicazioni}. }
}
\begin{document}
\maketitle
\begin{abstract}  
We prove that for every planar differential system with a period annulus there exists an involution $\sigma$ such that the system is 
$\sigma$-symmetric. We also prove that for  for every planar differential system with a period annulus there exist infinitely many involutions $\sigma$ such that the system is $\sigma$-reversible. 

 \end{abstract}



\section{Introduction}

Let us consider a planar differential system
\begin{equation} \label{sys}
\dot z = V(z),
\end{equation}  
where $V=(P(x,y),Q(x,y)) \in C^1 (\Omega,\R^2)$,  $\Omega \subset \R^2$ open and connected, $z=(x,y) \in \Omega$. We denote by $\phi(t,z)$ the local flow defined by (\ref{sys}) in $\Omega$.  
An isolated critical point $O$ of (\ref{sys}) is said to be a {\it center} if it has a neighbourhood covered with concentric nontrivial cycles. If $O$ is a center, we call period annulus $N_O$ the largest connected region covered with such concentric nontrivial cycles, such that $O \in \overline N_O$. More generally, we call {\it period annulus} every region (i. e. an open and connected subset of $\Omega$) $A$ covered with non-trivial cycles. Given a period annulus $A$, one can define on $A$ the so-called {\it period function } $T(z)$, by assigning to each point $z\in A$ its minimal positive period $T(z)$. The function $T(z)$ is a first integral of class $C^1$ of (\ref{sys}). Its properties has been widely studied for the information they provide about problems like the existence and multiplicity of solutions to some Dirichlet,  Neumann and mixed BVP's. The period function is also an essential tool in studying the stability properties of the non-trivial periodic solutions and some bifurcation problems \cite{BBLT, Ca, CLH, Co, CS, GV, LZ, L, LH, MMP,R,vW,Z}. 

In a recent paper \cite{S} the system (\ref{sys}) has been considered together a second one,
\begin{equation} \label{sysna}
\dot z = \alpha(z) V(z), \qquad \alpha \in C^0(\Omega,\R), \qquad \alpha(z) \neq 0,
\end{equation}  
having the same orbits, but not the same solutions. A sufficient condition for such systems to have the same period function on a common period annulus has been found. In order to prove such result,  (\ref{sys})  has been assumed to satisfy some kind of symmetry, i. e. point-symmetry and time-reversibility, both in an extended sense. 
In order to describe such properties, let $\sigma: \Omega \rightarrow \Omega$ be a homeomorhism such that $\sigma^2(z) = \sigma(\sigma(z)) =z$. Such maps are usually called {\it involutions}.
We say that the system (\ref{sys}) is {$\sigma$-reversible} if ther exists an involution $\sigma$ such that
$$
\sigma(\phi(t,z)) = \phi(-t,\sigma(z)),
$$
for all $t$ and $z$ such that the above equality makes sense.
Similarly, we say that the system (\ref{sys}) is {$\sigma$-symmetric} if 
$$
\sigma(\phi(t,z)) = \phi(t,\sigma(z)).
$$
Such flow symmetries can be revealed by examining the symmetry properties of the vector field. In fact,  if $\sigma \in C^1(\Omega,\Omega)$, then the $\sigma$-reversibility of (\ref{sys}) can be checked by verifying the relationship
$$
 V(\sigma(z)) = - J_\sigma(z) \cdot V(z) ,
$$
where $J_\sigma(z)$ is the Jacobian matrix of $\sigma$ at $z$. Moreover, the  $\sigma$-symmetry is equivalent to
$$
  V(\sigma(z)) = J_\sigma(z) \cdot V(z).
$$
Under a suitable condition on $\alpha$, (\ref{sys}) and (\ref{sysna}) have the same period function on a period annulus $A$ such that $\sigma(A) = A$ \cite{S}.

This leads naturally to ask whether a system with a period annulus $A$ be $\sigma$-reversible or $\sigma$-symmetric on $A$. Several results are known for the simplest case, where $\sigma$ is the mirror symmetry, in some cases giving procedures to find the reversibility line, i. e. the line of fixed points of $\sigma$ \cite{BL,BM,C,GM,vW}.  The interest in mirror symmetry is strictly related to the study of differential systems integrability. In this perspective, it is sufficient to prove orbital symmetry, considered in \cite{ACGG,GM}.

\medskip

In this paper we prove that every system with a  period annulus $A$ is both  $\sigma$-reversible and $\sigma$-symmetric on $A$. Both symmetries can be constructed by the flow. In section 1 we describe the procedure for the $\sigma$-symmetry, proving that  there exists exactly one involution $\sigma$ such that  every cycle in $A$ is $\sigma$-symmetric. If one drops the symmetry assumption on the single cycles, maitaining it on the period annulus, then there exist infinitely many distinct symmetries.

In section 2 we prove that every period annulus admits infinitely many involutions $\sigma$ such that the flow is  $\sigma$-reversible on $A$. In this case one can even arbitrarily choose the symmetry curve, i. e. the fixed points curve, always existing for such symmetries. 

Both in the symmetric case and in the reversible one we do not assume any non-degeneracy condition. 
Our results apply also to degenerate centers, which have unbounded period functions.  

\section{Existence of $\sigma$-symmetries}

We say that a regular curve $\delta\in C^1(I,\Omega)$, where $I$ is an open real interval, is a {\it section} of  (\ref{sys}), if its tangent  vector $\delta'(s)$ is transversal to $V(\delta(s))$ for all $s \in I$. If $A$ is an open period annulus of   (\ref{sys}), we say that $\delta$,  is a global section of  (\ref{sys}) on $A$ if $\delta$ meets every cycle of $A$. 

\begin{theorem} \label{teosymm}
Let $A$ be a period annulus of  (\ref{sys}). Then there exists a unique non-trivial involution $\sigma $ such that (\ref{sys}) is $\sigma$-symmetric and every cycle is $\sigma$-invariant.
\end{theorem}
{\it Proof.}   
Let us define $\sigma$ as follows,
$$
\sigma(z) := \phi \left( \frac{T(z)}{2},z \right).
$$
The map $\sigma$ is of class $C^1$ because it is the composition of  maps of class $C^1$. Moreover, one has
$$
\sigma(\sigma(z)) =  \phi \left( \frac{T(z)}{2},\phi \left( \frac{T(z)}{2},z \right) \right) = \phi \left( T(z),z \right) = z,
$$
hence $\sigma$ is an involution. As for the symmetry, one has
$$
\sigma(\phi(t,z)) = \phi \left( \frac{T(z)}{2},\phi(t,z) \right) =  \phi \left( t , \phi \left( \frac{T(z)}{2},z \right) \right) = 
\phi(t,\sigma(z)).
$$

Conversely, assume $\sigma$ to be an involution, (\ref{sys}) to be $\sigma$-symmetric and every cycle $\gamma$ to be $\sigma$-invariant.  Let us choose arbitrarily a cycle $\gamma^*$ in $A$ and a point $z^* \in \gamma^*$. By the invariance, the point $\sigma(z^*)$ belongs to $\gamma^*$, too, hence there exists $t^* \in (0,T(z^*))$ such that $\sigma(z^*)= \phi(t^*,z^*)$. Then, by the $\sigma$-symmetry, given any point $z=\phi(t,z^*) \in \gamma^*$, one has
$$
\sigma(z)= \sigma(\phi(t,z^*)) = \phi(t,\sigma(z^*)) =
 \phi(t+t^*,z^*) = \phi(t^*,\phi(t,z^*)) = \phi(t^*,z).
$$
Moreover, the map $\sigma$ is an involution, hence for all $z\in\gamma^*$
$$
z= \sigma^2(z) = \phi(2t^*,z),
$$
that implies
$$
 t^* = \frac {k\ T(z)}{2}, \qquad k \in \Z.
$$
The cycle $\gamma^*$ is $T(z)$-periodic, hence 
$$
\left\{  z =     \phi \left( \frac {   k\ T(z)}{2} , z\right), \  k \in \Z \right \} = \left\{ z, \phi\left( \frac{T(z)}{2} , z\right)\right \} .
$$
Hence, either $\sigma(z)=z$, or $\sigma(z) = \phi\left( \frac{T(z)}{2} , z\right)$. 
If $\sigma(z)=z$ for a single $z\in \gamma^*$, then by continuity, one has  $\sigma(z)=z$ for all $z \in\gamma^*$. Similarly, $\sigma(z)=z$ for all $z \in\gamma^*$ then by continuity, one has $\sigma(z)=z$  on all nearby cycles, hence $\sigma(z) = z$ on all of $A$, which is not compatible with the non-triviality of $\sigma$. Concluding, one has 
$$
\sigma(z) = \phi\left( \frac{T(z)}{2} , z \right)
$$
for all $z \in A$.
\hfill  $\clubsuit$ \\

\section{$\sigma$-reversibility}

\begin{theorem}  Let $A$ be a period annulus of  (\ref{sys}), and $\delta$ a global section of  (\ref{sys}) on $A$. Then there exists an involution $\sigma $ such that (\ref{sys}) is $\sigma$-reversible and $\delta$ is the fixed curve of $\sigma$.
\end{theorem}
{\it Proof.}   
As in the prof of theorem \ref{teosymm}, we denote by $T(z)$ the period of the cycle passing through the point $z \in A$.
We define a curve $\delta^*$ as follows, 
$$
\delta^* := \phi\left( \frac{T(z)}{2}, \delta \right).
$$
The curve $\delta^*$ is of class $C^1$. It is a global section of  (\ref{sys}) on $A$. For every $z \in A$, we set
$$
\gamma_+(z) = \left\{ \phi(t,z) : \ t\in \left( 0, \frac{T(z)}{2} \right) \right \}, \quad
\gamma_-(z) = \left\{ \phi(t,z) : \ t\in \left(  \frac{T(z)}{2} ,T(z) \right) \right \}.
$$
The we set
$$ 
A_+ =  \cup_{z \in \delta}\  \gamma_+(z)  , \qquad A_- = \cup_{z \in \delta}\  \gamma_-(z)  .
$$
One has $ A = A_-  \cup \delta \cup A_+ \cup   \delta^*$. 
Let us define a map $\tau :A \setminus \delta^* \rightarrow \R$ as follows,
$$
\tau (z) =\left\{
\matrix{ 
 \sup  \ \{ t\leq 0 : \phi(t,z) \in \delta  \},  &  z \in A^+ \setminus \delta^* \hfill  \cr 
 \inf \ \{ t\geq 0 : \phi(t,z) \in \delta  \}, &  z \in A^- \setminus \delta^*    . }
\right.
$$
By construction one has $\tau(z) = 0$ for every $z \in \delta$, and
$$
\tau(z) \in \left( -\frac{T(z)}{2}, \frac{T(z)}{2} \right), \qquad \forall z \in A \setminus \delta^*. 
$$ 
Moreover, if $z , \phi(t,z) \in A \setminus \delta^*$, then $\tau(\phi(t,z)) = \tau(z) - t$. In particular, $\tau(\phi(2\tau(z),z)) = - \tau(z) $.

We also define a map $\tau^* :A \setminus \delta \rightarrow \R$ as follows,
$$
\tau^* (z) =\left\{
\matrix{ 
 \sup  \ \{ t\leq 0 : \phi(t,z) \in \delta^*  \},  &  z \in A^- \setminus \delta \hfill  \cr 
 \inf \ \{ t\geq 0 : \phi(t,z) \in \delta^*  \}, &  z \in A^+ \setminus \delta    . }
\right.
$$
The map $\tau^*$ has similar properties to those of $\tau(z)$.
By using the flow's group property, it is easy to show that 
$$\left\{
\matrix{ 
  \tau^*(z)= \tau(z) + \frac{T(z)}2 ,  &  z \in A^+  \hfill  \cr 
  \tau^*(z)= \tau(z) - \frac{T(z)}2 , &  z \in A^-      . }
\right.
$$
Now we define a map $\sigma: A \rightarrow A$ as follows
\begin{equation} \label{Asigma}
\sigma(z) =
\left\{
\matrix{ 
\phi(2\tau(z),z) ,  &  z \in A \setminus \delta^* \hfill  \cr 
\phi(2\tau^*(z),z)  , &  z \in A \setminus \delta    . }
\right.
\end{equation}  
Such a definition is well-posed, since for every $z \not \in \delta \cup \delta^*$ one has
\begin{equation} \label{tautau}
\phi(2\tau^*(z),z) = \phi(2\tau(z),z).
\end{equation}  
In fact, if $z \in A_+$, one has
$$
\phi(2\tau^*(z),z) = \phi(2\tau(z) +  T(z) ,z)  = \phi(2\tau(z)  ,z),
$$
and for $z\in A_-$ one has
$$
\phi(2\tau^*(z),z) = \phi(2\tau(z) -  T(z) ,z)  = \phi(2\tau(z)  ,z).
$$
Similarly one proves that the equality holds also on $\delta \cup \delta^*$.

We claim that $\sigma$ is of class $C^1$ on $A$. It is sufficient to prove it in a nieghbourhood of $\delta$, then extend it to all of $A$ by using the fact that the flow is a diffeomorphism. Let $z_0$ be an arbitrary point of $\delta$. By the rectification theorem \cite{Ch}, there exists a local diffeomorphism $\Lambda: U_0 \rightarrow \R^2$,  defined in a suitable neighbourhood  $U_0 \subset A$ of $z_0$, such that $\Lambda$ takes (\ref{sys}) into the system
$$
\left\{
\matrix{ 
\dot x = 1   \hfill \cr 
\dot y = 0    .
 }  \right. 
$$
Possibly choosing a smaller $U_0$, we may reduce to the case $U_0 \cap \delta^* = \emptyset$. The curve $\delta$ is taken into a curve $\delta_0$, transversal to the above vector field, hence representable as the graph of a function $x = \alpha(y)$. The new transformation 
$$
\left\{
\matrix{ 
u = x - \alpha(y)    \cr 
v = y \hfill   .
 }  \right. 
$$
does not change the vector field,
$$
\left\{
\matrix{ 
\dot u = 1   \hfill \cr 
\dot v = 0    ,
 }  \right. 
$$ 
and takes $\delta_0$ into the line $u=0$. For the sake of simplicity, let us call again $\Lambda$ the composition of such two transformations, $\phi_0$ the local flow defined by the last vector field, $\phi_0(t,(u,v)) =( u+t,v)$, and $\delta_0$ the line $u=0$. Considering only $A \cap U_0$ and $\Lambda(A \cap U_0)$, one has 
$$
\tau(z) =  \sup_{ z \in A^+  \cap U_0 }  \ \{ t\leq 0 : \phi(t,z) \in \delta  \}   =  \sup_{\Lambda(z) \in   \Lambda(A^+ \cap U_0)}  \ \{ t\leq 0 : \phi_0(t,\Lambda(z)) \in \delta_0  \},  \quad 
 $$
Since, setting $(u,v) = \Lambda(z)$,  one has $\phi_0(t,\Lambda(z)) = \phi_0(t,(u,v)) = (u+t,v)$, one can write
$$
\tau(z)  =  \sup  \ \{ t\leq 0 : \phi_0(t,\Lambda(z)) \in \delta_0  \} = - u. 
 $$ 
 Similarly, one proves that 
 $$
\tau(z)  = \inf_{ z \in A^-  \cap U_0 }  \ \{ t\geq 0 : \phi(t,z) \in \delta  \}   =  \inf_{\Lambda(z) \in   \Lambda(A^- \cap U_0)}  \ \{ t\geq 0 : \phi_0(t,\Lambda(z)) \in \delta_0  \}.
 $$
As above, this leads to
 $$
\tau(z)  =   \inf_{ z \in A^-  \cap U_0 }  \ \{ t\geq 0 : \phi_0(t,\Lambda(z)) \in \delta_0  \} = - u. 
 $$
 This shows that $\tau(z)$ is a $C^1$ function, since it is the composition of $\Lambda$ and the projection on the first component multiplied by $-1$. As a consequence,  $\sigma(z)$ is  $C^1$ on $A \cap U_0$. This argument does not depend on the particular point $z_0$ on $\delta$, hence it covers a neighbourhood $V_0$ of $\delta$. 
 
The regularity of $\sigma$ in $A \setminus \delta^*$ comes from the fact that the flow is a diffeomorphism. The regularity of $\sigma$ on $ \delta^*$ can be proved by working on $\tau^*$ just as done on $\tau$, and using the equality (\ref{tautau}).
 
Finally, we prove that $\sigma$ is an involution and that that the flow is $\sigma$-reversible.
One has
$$
\sigma^2(z) = \sigma\Big( \phi(2\tau(z),z)  \Big) = \phi\Big( 2\tau(\phi(2\tau(z),z),\phi(2\tau(z),z)) \Big) =
$$
$$
= \phi \Big( -2\tau(z),\phi(2\tau(z),z) \Big) = z.
$$
In order to prove that the flow is $\sigma$-reversible, we can write
$$
\sigma \Big( \phi(t,z) \Big) = \phi\Big( 2 \tau(\phi(t,z)),\phi(t,z) \Big) = \phi\Big( 2(\tau(z) - t),\phi(t,z) \Big) = 
$$
$$
= \phi\Big( 2\tau(z) - 2t, \phi(t,z) \Big) = \phi\Big( 2\tau(z) - t, z \Big) = \phi\Big( -t,\phi(2\tau(z),z)) \Big) = \phi\Big( -t,\sigma(z) \Big).
$$
\hfill  $\clubsuit$ \\

Another way to prove the regularity of $\sigma$ consists in observing that, if for $w= (u,v) \in  \Lambda(A  \cap U_0)$, one defines
$$
\sigma_0(u,v) = (-u,v).
$$
then, for $z \in  A  \cap U_0$, one has
$$
\sigma(z) = \Lambda^{-1}(\sigma_0(\Lambda(z))).
$$

\enddocument